\documentclass[12pt]{article}
\usepackage{amsfonts,amssymb,amsmath,amsthm}

\theoremstyle{plain}
\newtheorem{question}{Question}
\newtheorem{lemma}{Lemma}[section]
\newtheorem{cor}[lemma]{Corollary}
\newtheorem{thm}[lemma]{Theorem}
\newtheorem{prop}[lemma]{Proposition}

\theoremstyle{definition}
\newtheorem{definition}[lemma]{Definition}
\newtheorem{example}[lemma]{Example}

\theoremstyle{remark}

\newtheorem*{remark}{Remark}

\newcommand{\ann}{{\rm ann}}
\newcommand{\Bezout}{B\'{e}zout }

\begin{document}

\title{A Characterization of Certain Morphic Trivial Extensions}
\author{Alexander J. Diesl, Thomas J. Dorsey and Warren Wm. McGovern}
\maketitle
\bigskip

\begin{abstract}
Given a ring $R$, we study the bimodules $M$ for which the trivial extension $R\propto M$ is morphic.  We obtain a complete characterization in the case where $R$ is left perfect, and we prove that $R\propto Q/R$ is morphic when $R$ is a commutative reduced ring with classical ring of quotients $Q$.  We also extend some known results concerning the connection between morphic rings and unit regular rings.

\end{abstract}

\section{Introduction}

A well-known result of G. Ehrlich \cite[Theorem 1]{ehr:unit} states that an endomorphism $\varphi$ of a module $M$ is unit regular if and only if $\varphi$ is von Neumann regular and $M/{\rm im}\ \varphi\cong{\rm ker}\ \varphi$. It is the latter condition that motivates the definition of a left morphic ring.   Following \cite{nisc:morp}, an element $a$ in a ring $R$ is called {\em left morphic} if $R/Ra\cong\ann_l(a)$ as left $R$-modules, and a ring $R$ is called left morphic if every element of $R$ is left morphic.  A ring that is both left and right morphic is called a {\em morphic ring}.  Ehrlich's theorem then implies that the class of left morphic rings includes all unit regular rings.

It is proved in \cite[Lemma 1]{nisc:morp} that an element $a$ in a ring $R$ is left morphic if and only if there is an element $b\in R$ such that $\ann_l(a)=Rb$ and $\ann_l(b)=Ra$.  Further generalizing, a {\em left quasi-morphic} element $a$ of a ring $R$ is defined in \cite{cn:qmr} to be one for which there are elements $b,c\in R$ such that $\ann_l(a)=Rb$ and $\ann_l(c)=Ra$.  A ring in which every element is left quasi-morphic is called a left quasi-morphic ring, and a ring which is both left and right quasi-morphic is called a {\em quasi-morphic ring}.  It is immediate that the class of left quasi-morphic rings contains all left morphic rings as well as all von Neumann regular rings.

The behavior of the morphic and quasi-morphic conditions in rings which satisfy certain chain conditions is studied in \cite{nisc:morp}, \cite{nisc:prin} and \cite{cn:qmr}.  A major consequence of these investigations is the characterization of the quasi-morphic rings that satisfy either the ascending or descending chain condition on principal left ideals as precisely the artinian rings in which every one-sided ideal is principal, the latter class having been described by Jacobson in \cite[Section 15]{jac:ring}.

If $R$ is a ring and ${}_RM_R$ is an ($R,R$)-bimodule, then the trivial extension $R\propto M$ is defined to be the set of all pairs $(r, m)$ such that $r\in R$ and $m\in M$.  Addition is defined componentwise, and multiplication is defined according to the rule $(r, m)(s, n)=(rs, rn+ms)$.  Given an endomorphism $\sigma$ of the ring $R$, we define the bimodule $R(\sigma )$ by ${}_RR(\sigma )={}_RR$ with right $R$-multiplication $s\cdot r=s\sigma (r)$ for every $s\in R(\sigma )$ and every $r\in R$.  Note that $R\propto R(\sigma )$ is isomorphic to $R[t;\sigma ]/(t^2)$ where $R[t;\sigma ]$ is a skew (left) polynomial ring over $R$.  When convenient, this isomorphism will hereafter be understood without further mention.

Many results on morphic trivial extensions are introduced in \cite{cz:ext} and further extended in \cite{lz:ext}.  In this article, we continue to extend the known results, focusing in particular on the cases where the base ring is either a unit regular ring, a left perfect ring or a commutative reduced ring.  Although our main results concern morphic trivial extensions, we shall often extend our results to the quasi-morphic case when it will not take us too far afield to do so.

All rings are assumed to be associative and unital.  The Jacobson radical and group of units of a ring $R$ will be denoted by $J(R)$ and $U(R)$, respectively.  Anytime that $R\propto M$ is written, it will be understood that ${}_RM_R$ is an ($R,R$)-bimodule.  A {\em principal ring} will refer to a ring which is both a principal left ideal ring and a principal right ideal ring.

\section{Preliminary Results}

The following fundamental result about left quasi-morphic rings is indispensable to our investigation.

\begin{prop}\label{prop:lqmbez}{\rm \cite[Corollary 7]{cnw:lqm}}
A left quasi-morphic ring is left B\'{e}zout.
\end{prop}

In particular, any left morphic ring is left B\'{e}zout.

Collecting \cite[Corollary 16]{nisc:prin} and \cite[Theorem 19]{cn:qmr} together with the classical structure theorem for artinian principal rings \cite[Section 15]{jac:ring}, we obtain the following result.

\begin{prop}\label{prop:quasi}{\rm (Cf. \cite[Corollary 16]{nisc:prin}, \cite[Theorem 19]{cn:qmr})}
The following are equivalent for a ring $R$.
\begin{enumerate}
\item[{\rm (1)}] $R$ is quasi-morphic and has either ACC or DCC on principal left ideals.
\item[{\rm (2)}] $R$ is an artinian principal ring.
\item[{\rm (3)}] $R$ is an artinian morphic ring.
\item[{\rm (4)}] $R$ is a morphic principal ring.
\item[{\rm (5)}] $R$ is a finite direct product of matrix rings over artinian principal local rings.
\end{enumerate}
\end{prop}

The following example shows that Proposition~\ref{prop:quasi} cannot be extended to rings which are left quasi-morphic but not right quasi-morphic.

\begin{example}{\rm (\cite[p. 70]{bjork:ex},\cite[Example 4]{cnw:lqm})}\label{ex:bjo}
Let $p$ be a prime, and let $F=\mathbb{F}_p(x)$.  Define $\sigma\colon F\to F$ by $\sigma(a)=a^p$, and let $R=F[t;\sigma ]/(t^2)$.  Then $R$ is a left morphic artinian ring which is not right quasi-morphic.
\end{example}

We now proceed with some general results about morphic trivial extensions.

\begin{prop}\label{prop:bez}
Let $R$ be a ring, and let $M$ be a bimodule.  If $S=R\propto M$ is left morphic, then $R$ is a left B\'{e}zout ring and ${}_RM$ is a B\'{e}zout {\rm (}left{\rm )} module.
\end{prop}

\begin{proof}
Since $S$ is left morphic, it is left B\'{e}zout by Proposition~\ref{prop:lqmbez}.  The ring $R$, being a quotient of $S$, must therefore also be left B\'{e}zout.  Consider two elements $m,n\in M$.  Since $S$ is left B\'{e}zout, there is an element $(a,x)\in S$ such that $S(0,m)+S(0,n)=S(a,x)$.  This implies immediately that $a=0$ and therefore that $Rm+Rn=Rx$.  Thus ${}_RM$ is B\'{e}zout.
\end{proof}

The following is a consequence of the proof of \cite[Proposition 27]{nisc:morp}.  We include the proof here for completeness.

\begin{prop}\label{prop:rl}
Suppose that $R$ is a ring and that $a\in R$ is a morphic element.  If  $b\in R$ is any element such that $Ra=\ann_l(b)$ and $Rb=\ann_l(a)$, then $aR = \ann_r(b)$ and $bR = \ann_r(a)$.
\end{prop}

\begin{proof}
Since $a$ is right morphic, $aR=\ann_r(c)$ for some $c\in R$.  Thus $\ann_r(\ann_l(a))= \ann_r(\ann_l(\ann_r(c))) = \ann_r(c)=aR$.  Then
 $\ann_r(b)=\ann_r(Rb)=\ann_r(\ann_l(a))=aR$.  A similar argument shows that
$\ann_r(a)=\ann_r(Ra)=\ann_r(\ann_l(b))=bR$.
\end{proof}

We now turn our attention to the interplay between elements of $R$ and elements of $M$ when $R\propto M$ is morphic.  Annihilators will play a central role.  We begin with a general result.

\begin{lemma}\label{lem:ann}
Let $R$ be a ring, let $M$ be a bimodule and let $S=R\propto M$.
\begin{enumerate}
\item[{\rm (A)}]The following are equivalent for elements $a\in R$ and $m,n\in M$:
\begin{enumerate}
\item[{\rm (1)}] $\ann_l^S(0,m)=S(a,n)$.
\item[{\rm (2)}] $\ann_l^R(m)=Ra$ and $\ann_l^R(a)n+Ma=M$.
\end{enumerate}
\item[{\rm (B)}]The following are equivalent for elements $a\in R$ and $m,n\in M$:
\begin{enumerate}
\item[{\rm (1)}] $\ann_l^S(a,n)=S(0,m)$.
\item[{\rm (2)}] $\ann_l^M(a)=Rm$, $\ann_l^R(a)n\cap Ma=0$ and $\ann_l^R(a)\cap\ann_l^R(n)=0$.
\end{enumerate}
\end{enumerate}
\end{lemma}

\begin{proof}
(A):  Suppose first that (1) holds.  Since $(a,n)(0,m)=0$, it is apparent that $Ra\subseteq\ann_l^R(m)$.  On the other hand, suppose that $b\in\ann_l^R(m)$.  Then $(b,0)\in\ann_l^S(0,m)$, which means that $(b,0)\in S(a,n)$, implying that $b\in Ra$.  Thus $\ann_l^R(m)=Ra$.  Suppose now that $x\in M$ is any element.  Since $(0,x)(0,m)=0$, there is an element $(c,y)\in S$ such that $(0,x)=(c,y)(a,n)=(ca, cn+ya)$.  Since $c\in\ann_l^R(a)$ and $y\in M$,  $x\in\ann_l^R(a)n+Ma$.  This shows that $\ann_l^R(a)n+Ma=M$.

On the other hand, suppose that (2) holds.  It is immediate that $S(a,n)\subseteq\ann_l^S(0,m)$.  Suppose that $(b,x)(0,m)=0$.  Since $bm=0$, we can write $b=ra$ for some $r\in R$, and since $M=\ann_l^R(a)n+Ma$, we can write $x-rn=cn+ya$ for some $c\in\ann_l^R(a)$ and some $y\in M$.  Then $(c+r,y)(a,n)=(ca+ra,cn+rn+ya)=(b,x)$.  Thus $\ann_l^S(0,m)=S(a,n)$.

(B):  Suppose that (1) holds.  It is immediate that $Rm\subseteq\ann_l^M(a)$.  Let $x\in\ann_l^M(a)$.  Then $(0,x)(a,n)=0$, which implies that $(0,x)\in S(0,m)$ and therefore that $x\in Rm$.  Thus $\ann_l^M(a)=Rm$.  Suppose now that $y\in\ann_l^R(a)n\cap Ma$.  Then $cn=y=za$ for some $c\in\ann_l^R(a)$ and some $z\in M$.  Then $(c,-z)(a,n)=(ca, cn-za)=0$, which implies that $(c,-z)\in S(0,m)$.  Thus $c=0$, which means that $y=0$.  We have then shown that $\ann_l^R(a)n\cap Ma=0$.  Finally, if $r\in\ann_l^R(a)\cap\ann_l^R(n)=0$, then $(r,0)(a,n)=(ra,rn)=0$.  Thus $(r,0)\in S(0,m)$, which implies that $r=0$.

On the other hand, suppose that (2) holds.  It is again clear that $S(0,m)\subseteq\ann_l^S(a,n)$.  Suppose that $(b,x)(a,n)=0$.  Then $ba=0$ and $bn+xa=0$.  But then $bn=-xa\in\ann_l^R(a)n\cap Ma=0$.  Thus $b\in\ann_l^R(a)\cap\ann_l^R(n)=0$ and $x\in\ann_l^M(a)=Rm$.  But then $(b,x)=(0,rm)=(r,0)(0,m)\in S(0,m)$ for some $r\in R$.  Thus $\ann_l^S(a,n)=S(0,m)$.
\end{proof}

\begin{cor}\label{lem:fund}Let $R$ be a ring and $M$ be a bimodule.  If $S=R\propto M$ is left morphic, then for every $m\in M$ there exists an element $a\in R$ such that $\ann_l^R(m)=Ra$ and $\ann_l^M(a)=Rm$.
\end{cor}

\begin{proof} 
Let $m\in M$, and write $S=R\propto M$.  Since $S$ is left morphic, there is an element $(a,n)\in S$ such that $\ann_l^S(0,m)=S(a,n)$ and $\ann_l^S(a,n)=S(0,m)$.  By Lemma~\ref{lem:ann},  $\ann_l^R(m)=Ra$ and $\ann_l^M(a)=Rm$.
\end{proof}

Combining Lemma~\ref{lem:ann} with its right-hand analog, we obtain the following result.

\begin{cor}\label{cor:fund}
Suppose that $S=R\propto M$ is a morphic ring.  For every $m\in M$ there is an $a\in R$ such that
$$\ann_l^R(m)=Ra\ \ \ \ \ \ \ \ \ \ \ \ann_l^M(a)=Rm$$
$$\ann_r^R(m)=aR\ \ \ \ \ \ \ \ \ \ \ \ann_r^M(a)=mR.$$
\end{cor}

\begin{proof}
Let $m\in M$.  Following the proof of Corollary~\ref{lem:fund}, there is an element $(a,n)\in S$ such that $\ann_l^S(0,m)=S(a,n)$ and $\ann_l^S(a,n)=S(0,m)$.  Since $S$ is also right morphic by hypothesis, Proposition~\ref{prop:rl} guarantees that $\ann_r^S(0,m)=S(a,n)$ and $\ann_r^S(a,n)=S(0,m)$ as well.  Using the right-hand analog of Lemma~\ref{lem:ann}, we obtain the desired equalities.
\end{proof}

Corollary~\ref{cor:fund} allows us to establish our first structure theorem.  For any B\'{e}zout ring $R$, let $\mathcal{L}(R)$ denote the lattice of principal left ideals of $R$, and for any B\'{e}zout right $R$-module $M$, let $\mathcal{R}(M)$ denote the lattice of cyclic right $R$-modules.  Corollary~\ref{lem:fund} shows that, if $R\propto M$ is left morphic, then we can define a function
$$\mathcal{F}\colon\mathcal{R}(M)\to\mathcal{L}(R)$$
by the rule $\mathcal{F}(N)=\ann_l^R(N)$.  This result and more are contained in the following theorem.

\begin{thm}\label{thm:str}
Suppose that $R$ is a ring and that $M$ is a bimodule such that $S=R\propto M$ is a morphic ring.
Then the function $\mathcal{F}$ defined above is an inclusion-reversing injective map.  If $N_R\in\mathcal{R}(M)$ is a sub-bimodule of $M$, then $\mathcal{F}(N)$ is an ideal of $R$.
\end{thm}

\begin{proof}
Suppose that $N_R$ is a cyclic right submodule of $M$.  Then $N_R$ has the form $mR$ for some $m\in M$.  By Corollary~\ref{lem:fund} there is an element $a\in R$ such that $\ann_l^R(mR)=\ann_l^R(m)=Ra$.  This
guarantees that the map $\mathcal{F}$ is well-defined.  It is straightforward to show that $\mathcal{F}$ is inclusion-reversing.

In order to prove that $\mathcal{F}$ is injective, suppose that there are elements $m,n\in M$ such that $\ann_l^R(mR)=\ann_l^R(nR)$.  Corollary~\ref{cor:fund} provides elements $a,b\in R$ such that $\ann_l^R(m)=Ra$, $\ann_r^M(a)=mR$, $\ann_l^R(n)=Rb$ and $\ann_r^M(b)=nR$.  If $\mathcal{F}(mR)=\mathcal{F}(nR)$, then $Ra=Rb$.  Applying $\ann_r^M$ to both sides, we see that $mR=nR$, which establishes the injectivity.

Finally, suppose that $mR\in\mathcal{R}(M)$ is a bimodule, and suppose that $\mathcal{F}(mR)=Ra$.  Since $mR$ is a bimodule, $rm\in mR$ for every $r\in R$.  Therefore, $(ar)m=a(rm)\in a(mR)=0$ for every $r\in R$.  Thus $aR\subseteq\ann_l^R(m)=Ra$, demonstrating that $Ra$ is an ideal.
\end{proof}

Theorem~\ref{thm:str} illustrates how the left and right submodule lattices of $M$ are influenced by the right and left ideal lattices of $R$ when $R\propto M$ is a morphic ring.  In what follows, we will exploit this relationship to illuminate the structure of $M$ in the cases where $R$ satisfies a chain condition or when $M$ is, itself, cyclic.

In order to state the next result, we need the following construction.  Let $I$ be an ideal of a ring $R$, and denote the quotient $R/I$ by $\overline{R}$.  For any ring endomorphism $\sigma\colon\overline{R}\to\overline{R}$, the ($\overline{R},\overline{R}$)-bimodule $\overline{R}(\sigma )$ can be viewed as an ($R,R$)-bimodule in the natural way.

\begin{thm}\label{thm:leftcyc}
Suppose that $R$ is a ring and that $M$ is a bimodule such that $R\propto M$ is left morphic.  Suppose further that $M=Rx$ for some $x\in M$ and that $\ann_l^R(x)$ is an ideal of $R$.  Then, writing $\overline{R}=R/\ann_l^R(x)$,  $M$ is isomorphic {\rm (}as a bimodule{\rm )} to $\overline{R}(\sigma )$ for some endomorphism $\sigma$ of $\overline{R}$.
\end{thm}

\begin{proof}
By Corollary~\ref{lem:fund}, there is an element $a\in R$ such that $\ann_l^R(x)=Ra$ (which is an ideal by hypothesis) and $\ann_l^M(a)=Rx$.  We now construct the endomorphism $\sigma$.  Since $M=Rx$, for every $r\in R$ there is an element $s\in R$ such that $xr=sx$.  Define a map $\varphi\colon R\to\overline{R}$ by $\varphi (r)=\overline{s}$ if $xr=sx$.  If $s$ and $s'$ are two elements of $R$ such that $sx=s'x$, then $(s-s')x=0$.  Since $Ra=\ann_l^R(x)$, $s-s'\in Ra$ and therefore $\overline{s}=\overline{s'}$.  Thus $\varphi$ is well-defined.  We now show that $\varphi$ is a ring homomorphism.  It is clear that $\varphi (1)=\overline{1}$.  If $r,r'\in R$ such that $\varphi (r)=\overline{s}$ and $\varphi (r')=\overline{s'}$, then $x(r+r')=xr+xr'=sx+s'x=(s+s')x$ and $x(rr')=(xr)r'=(sx)r'=s(xr')=s(s'x)=(ss')x$.  Therefore $\varphi (r+r')=\overline{s}+\overline{s'}=\varphi (r)+\varphi (r')$ and $\varphi(rr')=\overline{s}\overline{s'}=\varphi (r)\varphi (r')$.  Finally, we verify that ${\rm ker}\ \varphi\subseteq Ra$.  An element $r\in R$ is in ${\rm ker}\ \varphi$ if and only if $xr=sx$ for some $s\in Ra$.  Since $(Ra)x=0$, ${\rm ker}\ \varphi =\ann_r^R(x)$.  Since $xR\subseteq Rx=M$, $x(Ra)=(xR)a\subseteq Rxa=0$.  Thus $Ra\subseteq\ann_r^R(x)={\rm ker}\ \varphi$.  The map $\sigma\colon\overline{R}\to\overline{R}$ defined by $\sigma (\overline{r})=\varphi (r)$ is then the desired endomorphism.

It remains to be shown that $M$ is isomorphic, as a bimodule, to $\overline{R}(\sigma )$.  Fixing the element $x\in M$ such that $M=Rx$, define a map $\psi\colon\overline{R}(\sigma )\to M$ by $\psi (\overline{s})=sx$.  As above, we see that $\psi$ is a well-defined bijection.  It is also easy to see that $\psi$ is additive.  It remains to be shown that $\psi$ is both left and right $R$-linear.  Let $a\in R$ and $\overline{s}\in\overline{R}(\sigma )$.  Then
$$\psi (a\cdot\overline{s})=\psi (\overline{a}\overline{s})=(as)x=a(sx)=a\psi (\overline{s})$$
and
$$\psi (\overline{s}\cdot a)=\psi (\overline{s}\sigma (\overline{a}))=(s\varphi (a))x= s(\varphi (a)x)=s(xa)=(sx)a=\psi (\overline{s})a.$$
\end{proof}

If $R\propto M$ is morphic, we can obtain a stronger result.

\begin{thm}\label{thm:quot}
Suppose that $R$ is a ring and that $M$ is a bimodule such that $R\propto M$ is morphic.  If $M_R$ and ${}_RM$ are cyclic, then there is an ideal $I$ of $R$ and an automorphism $\sigma$ of $\overline{R}=R/I$ such that $M\cong\overline{R}(\sigma )$.
\end{thm}

\begin{proof}
Since ${}_RM$ and $M_R$ are cyclic, there are elements $x,y\in M$ such that ${}_RM=Rx$ and $M_R=yR$.  Corollary~\ref{cor:fund} provides elements $a,b\in R$ such that
$$\ann_l^R(x)=Ra\qquad\ann_l^M(a)=Rx$$
$$\ann_r^R(x)=aR\qquad\ann_r^M(a)=xR$$
and
$$\ann_l^R(y)=Rb\qquad\ann_l^M(b)=Ry$$
$$\ann_r^R(y)=bR\qquad\ann_r^M(b)=yR.$$

The proof will be a strengthening of the proof of Theorem~\ref{thm:leftcyc}.  We first claim that $Rx=xR$ (and similarly that $Ry=yR$).  Since $M$ is a bimodule, it is immediate that $xR\subseteq Rx$ and $Ry\subseteq yR$.  Since $Rx$ and $yR$ are bimodules, Theorem~\ref{thm:str} (and its analog on the right) guarantees that $aR$ and $Rb$ are ideals of $R$.  This implies that $Ra\subseteq aR$ and $bR\subseteq Rb$.  Since ${}_RM=Rx$, $Ry\subseteq Rx$.  Taking right annihilators in $R$, we see that $aR=\ann_r^R(x)\subseteq\ann_r^R(y)=bR$.  A similar argument shows that the inclusion $xR\subseteq yR$ implies that $Rb\subseteq Ra$.  Combining these inclusions yields
$$Ra\subseteq aR\subseteq bR\subseteq Rb\subseteq Ra.$$
Thus $Ra=aR=Rb=bR$; similarly $Rx=xR=Ry=yR$.

Let $\overline{R}$ denote $R/Ra$.  Define the map $\varphi\colon R\to\overline{R}$ as in the proof of Theorem~\ref{thm:leftcyc}.  It remains only show that $\varphi$ is surjective and that ${\rm ker}\ \varphi =Ra$.  Since $M=xR$, we see immediately that $\varphi$ is surjective.  As in the proof of Theorem~\ref{thm:leftcyc}, ${\rm ker}\ \varphi =\ann_r^R(x)$; since $\ann_l^R(x)=aR=Ra$, the map $\sigma\colon\overline{R}\to\overline{R}$ defined by $\sigma (\overline{r})=\varphi (r)$ is the desired automorphism.
\end{proof}

We remark that, by Proposition~\ref{prop:bez}, the conclusions of Theorem~\ref{thm:leftcyc} and Theorem~\ref{thm:quot} hold more generally in the case where the relevant modules are finitely generated.

\begin{cor}\label{cor:tor}
Suppose that $R$ is a ring and that $M$ is a bimodule such that $R\propto M$ is left morphic.  If there is an element $x\in M$ such that $\ann_l^R(x)=0$, then $R\propto M\cong R[t;\sigma]/(t^2)$ for some endomorphism $\sigma$ of $R$.  In particular, if $R[t;\sigma ]/(t^2)$ is morphic, then $\sigma$ must be an automorphism.
\end{cor}

\begin{proof}
By Corollary~\ref{lem:fund}, since $\ann_l^R(x)=0$, we must have $Rx=\ann_l^M(0)={}_RM$.  By Theorem~\ref{thm:leftcyc}, $M\cong R(\sigma )$ for some endomorphism $\sigma$ of $R$.  Thus $R\propto M\cong
R\propto R(\sigma )\cong R[t;\sigma ]/(t^2)$.  If $R\propto M$ is morphic, then Corollary~\ref{cor:fund} implies that $xR=\ann_r^M(0)=M_R$, and Theorem~\ref{thm:quot} guarantees that $\sigma$ is an automorphism.
\end{proof}

The next result is a generalization of \cite[Theorem 19]{cz:ext}.

\begin{thm}\label{thm:gencz}
Suppose that $R$ is a ring and that $\sigma$ is an endomorphism of $R$.  
\begin{enumerate}
\item[(1)] If $(a,0)$ is left morphic in $R\propto R(\sigma )$ and $\sigma$ is an automorphism, then $a$ is left morphic in $R$.
\item[(2)] If $(a,0)$ is morphic in $R\propto R(\sigma )$, then $a$ is morphic in $R$.
\end{enumerate}
\end{thm}

\begin{proof}
Let $S$ denote $R\propto R(\sigma )$.  If $(a,0)$ is left morphic in $S$, then there is an element $(b,c)\in S$ such that $\ann_l^S(a,0)=S(b,c)$ and $\ann_l^S(b,c)=S(a,0)$.  Since $(a,0)(b,c)=(0,0)=(b,c)(a,0)$, it is immediately clear that $ab=0=ba$ and $ac=0=c\sigma (a)$.  Thus, $Rb\subseteq\ann_l^R(a)$ and $Ra\subseteq\ann_l^R(b)$.  Further, if $xa=0$, then $(x,0)\in\ann_l^S(a,0)=S(b,c)$, implying that $x\in Rb$ and therefore that $Rb=\ann_l^R(a)$.

If $\sigma$ is an automorphism, then $R\sigma (b)=\ann_l^R(\sigma (a))$.  Since $c\sigma (a)=0$,  $c\in R\sigma (b)$.  Thus $c=r\sigma (b)$ for some $r\in R$.  We then note that $(b,c)=(b, r\sigma (b))=(1,r)(b,0)$.  Since $(1,r)$ is a unit in $S$ with inverse $(1,-r)$, \cite[Lemma 3]{nisc:morp} implies that $\ann_l^S(b,0)=S(a,0)(1,-r)=S(a,-ar)$.  Thus if $y\in\ann_l^R(b)$, then $(y,0)\in\ann_l^S(b,0)=S(a,-ar)$, which implies that $y\in Ra$.  This shows that $Ra=\ann_l^R(b)$ and proves (1).

If, on the other hand, $(a,0)$ is morphic, then Proposition~\ref{prop:rl} implies that $\ann_r^S(a,0)=S(b,c)$ and $\ann_r^S(b,c)=S(a,0)$.  As above, $bR=\ann_r^R(a)$ and $aR\subseteq\ann_r^R(b)$. Since $ac=0$, $c\in\ann_r^R(a)=bR$.  Thus $c=by$ for some $y\in R$.  Thus $(b,c)=(b,0)(1,y)$; since $(1,y)$ is a unit in $S$, \cite[Lemma 3]{nisc:morp} implies that $\ann_l^S(b,0)=S(1,y)(a,0)=S(a,y\sigma (a))$.  If $z\in\ann_l^R(b)$, then $(z,0)\in\ann_l^S(b,0)=S(a,y\sigma (a))$, which implies that $z\in Ra$.  Therefore $Ra=\ann_l^R(b)$, and similarly, $aR=\ann_r^R(b)$.  Thus $a$ is morphic in $R$, which proves (2).
\end{proof}

\begin{cor}\label{cor:gencz}
If $R$ is a ring and $\sigma$ is an endomorphism of $R$ such that $R\propto R(\sigma )$ is morphic, then $R$ is morphic and $\sigma$ is an automorphism.
\end{cor}

\begin{proof}
Corollary~\ref{cor:tor} implies that $\sigma$ is an automorphism, and Theorem~\ref{thm:gencz} then shows that $R$ is morphic.
\end{proof}

Note that Corollary~\ref{cor:gencz} offers another proof that the ring in Example~\ref{ex:bjo} is not morphic since the endomorphism in that example is not an automorphism.

The next results are a direct consequence of Lemma~\ref{lem:ann}, and generalize \cite[Theorem 9]{lz:ext}.

\begin{cor}\label{cor:reg}
Let $R$ be a ring, and let $\sigma$ be an endomorphism of $R$.  For any $a\in R$, if $(0,a)$ is left morphic in $R\propto R(\sigma)$, then $a$ is von Neumann regular in $R$.
\end{cor}

\begin{proof}
Let $S=R\propto R(\sigma)$.  Since $(0,a)$ is left morphic in $S$, there is an element $(b,c)\in S$ such that $\ann_l^S(0,a)=S(b,c)$ and $\ann_l^S(b,c)=S(0,a)$.  Using Lemma~\ref{lem:ann} (and viewing all multiplications in $R$), we see that $\ann_l^R(b)c\oplus R\sigma (b)=R$.  Thus $\sigma (b)$ is von Neumann regular.  Appealing again to Lemma~\ref{lem:ann}, we see that $Ra=\ann_l^R(\sigma(b))$, which shows that $a$ is von Neumann regular.
\end{proof}

Although we will only need the next result in a particular case, the proof of the general case is no more difficult, and we include it here.

\begin{lemma}\label{lem:ureg}
Suppose that $R$ is a ring and that $M$ is a bimodule.  If $a$ is von Neumann regular in $R$ and $(a,0)$ is left morphic in $R\propto M$, then $a$ is unit regular.
\end{lemma}

\begin{proof}
Since $a$ is von Neumann regular, $a=axa$ for some $x\in R$.  Then $(a,0)(x,0)(a,0)=(a,0)$, so $(a,0)$ is regular in $R\propto M$.  Since $(a,0)$ is also left morphic in $R\propto M$, $(a,0)$ is unit regular.  Thus there is a unit $(u,m)\in R\propto M$ such that $(a,0)=(a,0)(u,m)(a,0)=(aua,0)$.  Since $u$ must be a unit in $R$, this shows that $a$ is unit regular.
\end{proof}

\begin{cor}\label{cor:ureg}
Let $R$ be a ring, and let $\sigma$ be an endomorphism of $R$.  If $R\propto R(\sigma)$ is left morphic, then $R$ is unit regular.
\end{cor}

\begin{proof}
By Corollary~\ref{cor:reg}, $R$ is von Neumann regular.  Lemma~\ref{lem:ureg} then shows that $R$ is unit regular.
\end{proof}

We note that Corollary~\ref{cor:ureg} answers the question following Theorem 9 in \cite{lz:ext}.

\begin{cor} Suppose that $R$ is a ring and that $\sigma$ is an endomorphism of $R$ that fixes all idempotents of $R$.  Then the following conditions are equivalent.
\begin{enumerate}
\item[{\rm (1)}] $R\propto R(\sigma )$ is left morphic.
\item[{\rm (2)}] $R$ is unit regular.
\end{enumerate}
\end{cor}

\begin{proof}
By Corollary~\ref{cor:ureg}, (1) implies (2).  The reverse implication holds by \cite[Theorem 2]{lz:ext}.
\end{proof}

Recall that a ring is called {\em strongly morphic} if $\mathbb{M}_n(R)$ is morphic for every positive integer $n$.

\begin{cor}
The following are equivalent for a ring $R$.
\begin{enumerate}
\item[{\rm (1)}] $R\propto R$ is left morphic.
\item[{\rm (2)}] The element $(0,a)$ is left morphic in $R\propto R$ for every $a\in R$.
\item[{\rm (3)}] $R$ is unit regular.
\item[{\rm (4)}] $R\propto R$ is morphic.
\item[{\rm (5)}] $R[x]/(x^n)$ is strongly morphic for all $n\ge 1$.
\end{enumerate}
\end{cor}

\begin{proof}
The equivalence of (3), (4) and (5) is proved in \cite[Theorem 9]{lz:ext}, and it is clear that (4)$\Rightarrow$(1)$\Rightarrow$(2).  In order to prove that (2) implies (3), let $a\in R$.  Since $(0,a)$ is left morphic in $R\propto R$, Lemma~\ref{lem:ann} shows that $a$ is left morphic in $R$, and Corollary~\ref{cor:reg} shows that $a$ is von Neumann regular in $R$.  Therefore, $a$ is unit regular in $R$.
\end{proof}

We close the section with some general results on idempotents.  Such results will be useful both in establishing a partial converse to Corollary~\ref{cor:ureg} as well as in characterizing morphic extensions in the next section.

\begin{prop}\label{prop:idem}
Suppose that $R$ is a ring and that $M$ is a bimodule.  If $e$ and $e'$ are idempotents in $R$, then $(e,0)$ and $(e',0)$ are idempotents in $R\propto M$, and $(e,0)(R\propto M)(e',0)= eRe'\propto eMe'$.
\end{prop}

\begin{proof}
The proof is a straightforward calculation and is omitted.
\end{proof}

\begin{prop}\label{prop:prod}
Suppose that $R$ is a ring and $M$ is a bimodule such that $R\propto M$ is left morphic.  If $e\in R$ is a central idempotent, then $em=me$ for every $m\in M$.
\end{prop}

\begin{proof}
By Proposition~\ref{prop:idem}, $(e,0)(R\propto M)(1-e,0)=eR(1-e)\propto eM(1-e)$.  Since $e$ is central, $eR(1-e)=0$.  Thus $(e,0)R\propto M(1-e,0)=0\propto eM(1-e)$ is contained in the Jacobson radical of $R\propto M$.  By \cite[Corollary 19]{nisc:morp}, $(e,0)(R\propto M)(1-e,0)=0=(1-e,0)(R\propto M)(e,0)$.  In particular, $em(1-e)=0=(1-e)me$ for every $m\in M$.  Thus $em=eme=me$.
\end{proof}

\begin{remark}
Identifying $R[t;\sigma ]/(t^2)$ with $R\propto R(\sigma )$, Proposition~\ref{prop:prod} implies that if $R$ is a ring and $\sigma\colon R\to R$ is an endomorphism such that $R[t;\sigma ]/(t^2)$ is left morphic, then $\sigma (e)=e$ for every central idempotent $e\in R$.
\end{remark}

It is shown in \cite[Theorem 2]{lz:ext} that the ring $R[t;\sigma ]/(t^2)$ is left morphic if $R$ is unit regular and $\sigma$ fixes every idempotent of $R$.  In light of the above remark, it is, indeed, necessary that $\sigma$ fix every {\em central} idempotent of $R$.  When $R$ is semisimple, the condition that $\sigma$ fixes every central idempotent of $R$ is also sufficient, as is shown in \cite[Corollary 16]{dor:morp} (see also Theorem~\ref{thm:con} and Corollary~\ref{cor:char} below).  In particular, if $R$ is simple artinian, then $R[t;\sigma ]/(t^2)$ is morphic for {\em any} automorphism $\sigma$ of $R$, although there are any automorphisms (e.g. conjugation by a noncentral matrix) that do not fix all idempotents of $R$.

\begin{cor}
If $R$ is Boolean and $\sigma$ is an endomorphism of $R$, then $R[t;\sigma ]/(t^2)$ is left morphic if and only if $\sigma$ is the identity.
\end{cor}

\begin{proof}
Every element of a Boolean ring is a central idempotent.
\end{proof}

\begin{cor}
If $R$ is a strongly regular ring and $\sigma$ is an endomorphism of $R$, then $R[t;\sigma ]/(t^2)$ is left morphic if and only if $\sigma (e)=e$ for every idempotent $e\in R$.
\end{cor}

\begin{proof}
The forward direction is a consequence of the remark following Proposition~\ref{prop:prod}; the reverse is \cite[Theorem 1]{cz:ext}.
\end{proof}

\section{Rings Satisfying Chain Conditions}

Proposition~\ref{prop:quasi} illustrates the behavior of the quasi-morphic condition under the application of chain conditions. The main goal of this section will be to use Proposition~\ref{prop:quasi} to completely characterize the quasi-morphic trivial extensions $R\propto M$ when $R$ is a left perfect ring.  We begin with some results to show that the lattices of left and right submodules of $M$ have a particularly nice form in this case.

\begin{lemma}\label{lemma:perf}
If $R$ is a left perfect ring and ${}_RM$ is a B\'{e}zout module, then ${}_RM$ is noetherian and cyclic.
\end{lemma}

\begin{proof}
Since $R$ is left perfect, a result of Jonah \cite[Main Theorem]{jon:perf} guarantees that ${}_RM$ satisfies the ascending chain condition on cyclic submodules.  By Proposition~\ref{prop:bez}, every finitely generated submodule of ${}_RM$ is cyclic, and ${}_RM$ thus satisfies the ascending chain condition on finitely generated submodules.  It is a straightforward exercise to show that this implies that ${}_RM$ is noetherian.  Since ${}_RM$ is B\'{e}zout, it must therefore be cyclic.
\end{proof}

\begin{cor}\label{cor:perf}
Suppose that $R$ is a perfect ring and that $M$ is a bimodule such that $R\propto M$ is quasi-morphic. Then  ${}_RM$ and $M_R$ are cyclic and have finite length.
\end{cor}

\begin{proof}
By Proposition~\ref{prop:bez} and Lemma~\ref{lemma:perf} (and their analogs on the right), both ${}_RM$ and $M_R$ are noetherian and cyclic.  Since ${}_RM$ and $M_R$ are B\'{e}zout, every submodule of either ${}_RM$ or $M_R$ is therefore cyclic.  Since $R$ is both left and right perfect, $M_R$ and ${}_RM$ satisfy the descending chain condition on cyclic submodules (see, for instance, \cite[p. 344]{lam:fc}).  Thus ${}_RM$ and $M_R$ are both artinian and therefore both have finite length.
\end{proof}

We will also use the following general fact.

\begin{lemma}\label{lemma:perfext}
If $R$ is a left perfect ring and $M$ is a bimodule, then $R\propto M$ is left perfect.
\end{lemma}

\begin{proof}
Since $J(R\propto M)=J(R)\propto M$, we see immediately that $R\propto M/J(R\propto M)\cong R/J(R)$ is semisimple.  We claim now that $J(R\propto M)$ is left T-nilpotent.  Suppose that $\lbrace (r_i,m_i)\rbrace$ is a sequence of elements in $J(R\propto M)$.  Since $r_i\in J(R)$ and $R$ is left perfect, there is an index $n$ for which $r_1r_2\dots r_n=0$.  Then $(r_1,m_1)(r_2,m_2)\dots (r_n,m_n)=(0,m)$ for some $m\in M$.  Using again the fact that $R$ is left perfect, there is a $k$ such that $r_{n+1}r_{n+2}\dots r_{n+k}=0$.  Then $(r_1,m_1)\dots (r_{n+k},m_{n+k})=(0,m)(r_{n+1},m_{n+1})\dots (r_{n+k},m_{n+k})=(0,mr_{n+1}\dots r_{n+k})=0$.
\end{proof}

We can now state our characterization.

\begin{thm}\label{thm:cyc}
Suppose that $R$ is a left perfect ring and that $M$ is a bimodule such that $R\propto M$ is quasi-morphic.  Then $R$ is an artinian principal ring, and $M\cong\overline{R}(\sigma )$ for some quotient $\overline{R}$ of $R$ and some automorphism $\sigma$ of $\overline{R}$.  Further, $R\propto M$ is an artinian principal ring.
\end{thm}

\begin{proof}
Since $R$ is left perfect, Lemma~\ref{lemma:perfext} implies that $R\propto M$ is left perfect.  Proposition~\ref{prop:quasi} then shows that $R\propto M$ must be an artinian principal, hence morphic, ring.  This shows that $R$, being a quotient of $R\propto M$, must be an artinian principal ring.

By Corollary~\ref{cor:perf}, $M_R$ and ${}_RM$ are cyclic.  Theorem~\ref{thm:quot} then shows that there must be a quotient $\overline{R}$ of $R$ and an automorphism $\sigma$ of $\overline{R}$ such that $M\cong\overline{R}(\sigma )$.
\end{proof}

\begin{remark}
Theorem~\ref{thm:cyc} shows that if $R$ is left perfect, then $R\propto M$ is morphic if and only if it is quasi-morphic.  We will therefore concern ourselves only with the morphic case.
\end{remark}

We are now in a position to completely characterize the morphic trivial extensions of an arbitrary left perfect ring.  The next result reduces the problem to the case of a ring with no nontrivial central idempotents.

\begin{lemma}\label{lemma:prod}
Suppose that $R$ is a ring and $M$ is a bimodule such that $R\propto M$ is left morphic.  If $R=R_1\times R_2$ is the direct product of two rings, then there are bimodules $M_1$ and $M_2$ such that $R\propto M\cong (R_1\propto M_1)\times (R_2\propto M_2)$.
\end{lemma}

\begin{proof}
Let $e=(1,0)$ and $f=(0,1)$ in $R_1\times R_2=R$, and let $M_1=eM$ and $M_2=fM$.  Applying Proposition~\ref{prop:prod}, we see that $M_1$ is an ($R_1,R_1$)-bimodule, that $M_2$ is an ($R_2,R_2$)-bimodule, and that the map that takes $(r,m)\in R\propto M$ to $((er,em),(fr,fm))$ is the desired isomorphism.
\end{proof}

\begin{thm}\label{thm:con}
Suppose that $R$ is a left perfect ring with no nontrivial central idempotents and that $M$ is a bimodule.  If $R$ is simple artinian, then $R\propto M$ is morphic if and only if either $M=0$ or $M\cong R(\sigma )$ for some automorphism $\sigma$ of $R$.  If $R$ is not simple artinian, then $R\propto M$ is morphic if and only if $R$ is an artinian principal ring and $M=0$.
\end{thm}

\begin{proof}
Suppose $R$ is simple artinian.  If $R\propto M$ is morphic, then Theorem~\ref{thm:cyc} implies that $M\cong\overline{R}(\sigma )$ for some quotient $\overline{R}$ of $R$ and some automorphism $\sigma$ of $\overline{R}$.  Since $R$ is simple, either $\overline{R}=0$ or $\overline{R}=R$.  On the other hand, if $M=0$ then $R\propto M\cong R$ is morphic, and if $M\cong R(\sigma )$ for some automorphism $\sigma$ of $R$, then $R\propto M$ is morphic by \cite[Corollary 16]{dor:morp}).

We now handle the case where $R$ is not simple artinian.  Let $S=R\propto M$. Suppose that $S$ is morphic but that $M$ is nonzero.  Since $R$ is a left perfect ring with no nontrivial central idempotents, then $S$ must be an artinian principal ring by Theorem~\ref{thm:cyc}.  Since $S$ also has no nontrivial central idempotents, $S$ must therefore be isomorphic to $\mathbb{M}_n(L)$ for some artinian local principal ring $L$ by Proposition~\ref{prop:quasi}.  Write $J=J(S)$.  It is clear that $J=J(R)\propto M$.  Since the only ideals of $L$ are powers of $J(L)$ (see \cite[Theorem 9]{nisc:morp}), the only ideals of $S$ are of the form $J^i$.  Thus $0\propto M=J^r$ for some $r$.  Since $R$ is not simple artinian, $J(R)\ne 0$, and it must therefore be the case that $r\ge 2$.  Note that $0\propto MJ(R)=J^rJ=J^{r+1}=JJ^r=0\propto J(R)M$.  Since $J^r\ne 0$, $J^{r+1}\ne J^r$.  Thus $J(R)M=MJ(R)\ne M$.  Then, $J^2=J(R)^2\propto (J(R)M+MJ(R))$ cannot contain $J^r=0\propto M$, a contradiction.  We therefore conclude that $M=0$.

For the converse, Proposition~\ref{prop:quasi} shows that any artinian principal ring is morphic.
\end{proof}

We now give the complete characterization.

\begin{cor}\label{cor:char}
Suppose that $R$ is a left perfect ring and that $M$ is a bimodule such that $R\propto M$ is morphic.  Then  $R=R_1\times\dots\times R_n$ where each $R_i$ is a matrix ring over an artinian principal local ring, and $M$ can be written as $M=M_1\oplus\dots\oplus M_n$, where each $M_i$ is a bimodule over $R_i$.  Further, $M_i=0$ or $M_i\cong R_i(\sigma_i)$ for some automorphism $\sigma_i$ of $R_i$ if $R_i$ is simple artinian and that $M_i=0$ otherwise.  Conversely, any such trivial extension is morphic.
\end{cor}

\begin{proof}This follows from Lemma~\ref{lemma:prod} and Theorem~\ref{thm:con}.
\end{proof}

Since any proper quotient of a PID is artinian, Corollary~\ref{cor:char} generalizes  \cite[Theorem 16]{lz:ext}.

\begin{cor}\label{cor:noeth}
Suppose that $R$ is a ring and $M$ is a bimodule such that $R\propto M$ is quasi-morphic. If $R$ satisfies the ascending chain condition on principal left ideals and ${}_RM$ is finitely generated, then $R$ is an artinian principal ring and $R\propto M$ is morphic.  In particular, $M$ is a bimodule of the type described in Corollary~\ref{cor:char}.
\end{cor}

\begin{proof}
Since $R\propto M$ is quasi-morphic, $R$ is B{\' e}zout by Proposition~\ref{prop:bez}.  Thus $R$ satisfies the ascending chain condition on finitely generated left ideals and is therefore left noetherian.  Since ${}_RM$ is finitely generated, it is also noetherian.  Therefore $R\propto M$ is left noetherian.  By Proposition~\ref{prop:quasi}, $R\propto M$ is an artinian principal ring.  The remaining conclusions follow from Theorem~\ref{thm:cyc} and Corollary~\ref{cor:char}.
\end{proof}

The next example shows that the condition that ${}_RM$ is finitely generated cannot be eliminated from the hypotheses of Corollary~\ref{cor:noeth}.

\begin{example}\cite[Theorem 14]{cz:ext} If $R=\mathbb{Z}$ and $M=\mathbb{Q}/\mathbb{Z}$, then $R\propto M$ is a morphic trivial extension in which $R$ is noetherian but not artinian, and $M$ is not finitely generated.
\end{example}

\begin{cor}\label{cor:strong}
Suppose that $R$ is a ring and $M$ is a bimodule such that $R\propto M$ is morphic.  If $R$ is left perfect or if $R$ satisfies the ascending chain condition on principal left ideals and ${}_RM$ is finitely generated, then $R\propto M$ is strongly morphic.
\end{cor}

\begin{proof}
In either case, $R\propto M$ is an artinian principal ring, in which case $\mathbb{M}_n(R\propto M)$ is also an artinian principal ring, and therefore a morphic ring.
\end{proof}

\section{Commutative Reduced Rings}

In this section, we investigate morphic trivial extensions $R\propto M$ where $R$ is a commutative reduced ring and ${}_RM_R$ is a bimodule such that $mr=rm$ for all $r\in R$ and all $m\in M$; this last condition is equivalent to saying that $R\propto M$ is a commutative ring.  Since it is proved in \cite[Corollary 4]{cn:qmr} that every commutative quasi-morphic ring is morphic, we lose no generality by focusing on morphic rings in this case.

It is shown in \cite[Theorem 14]{cz:ext} under these hypotheses on $M$ that $\mathbb{Z}\propto M$ is morphic if and only if $M\cong\mathbb{Q}/\mathbb{Z}$, and this result is generalized to arbitrary PIDs in \cite[Theorem 14]{lz:ext}.
In this section, we will use the results of Proposition~\ref{prop:bez} and specialize Lemma~\ref{lem:ann} and Theorem~\ref{thm:str} in order to extend some of these results to arbitrary commutative reduced rings.  Some results for arbitrary domains will also be proved.

We first establish a general fact which will be useful in what follows.

\begin{lemma}\label{lem:prin}
Let $R$ be a commutative ring, and let $M$ be a bimodule such that $R\propto M$ is morphic.  For every $a\in R$, $\ann^R(a)$ is a principal ideal.
\end{lemma}

\begin{proof}
Let $a\in R$, and write $S=R\propto M$.  Since $S$ is morphic, there is an element $(b,m)\in S$ such that $\ann^S(a,0)=S(b,m)$.  It is immediate that $Rb\subseteq\ann^R(a)$.  On the other hand, if $ca=0$, then $(c,0)\in\ann^R(a,0)=S(b,m)$.  Thus $\ann^R(a)=Rb$.
\end{proof}

\begin{definition}
A ring is called left Rickart (or left p.p.) if the left annihilator of every element is generated by an idempotent.  A commutative Rickart ring $R$ is often called a weak Baer ring.
\end{definition}

\begin{remark}
Note that a \Bezout ring is weak Baer ring if and only if it is semihereditary.
\end{remark}

\begin{lemma}\label{lem:baer}
Suppose that $R$ is a commutative reduced B\'{e}zout ring.  If the annihilator of every element is principally generated, then $R$ is weak Baer.
\end{lemma}

\begin{proof}
Suppose that $a\in R$.  Let $\ann (a)=bR$ and $\ann (b)=cR$.  Note that we then have $\ann(c)=bR$.  Since $R$ is B\'{e}zout, $bR+cR=dR$ for some $d\in R$.  If $x\in bR\cap cR$, then $x^2=0$.  Since $R$ is reduced, this means that $bR\cap cR=0$.  Thus $bR\oplus cR=dR$, and $d$ is a nonzerodivisor since $\ann (d)=\ann (bR+cR)=\ann (b)\cap\ann (c)=cR\cap bR=0$.

Therefore $dR$ is a free $R$-module, which shows that $cR$ is also projective.  Thus the exact sequence
$$0\to\ann (c)\to R\to cR\to 0$$
splits, which shows that $\ann (a)=\ann (c)$ is a direct summand of $R_R$.
\end{proof}

If $R$ is a commutative ring, we will let $Q(R)$ (or simply $Q$ if there is no risk of confusion) denote the classical ring of quotients of $R$.

\begin{lemma}\label{lem:comm1}
Let $R$ be a commutative ring with classical ring of quotients $Q$.  If $a\in R$ is not a zerodivisor, then $\ann^R(\overline{\frac{1}{a}})=Ra$ and $\ann^{Q/R}(a)=R\overline{\frac{1}{a}}$.
\end{lemma}

\begin{proof}
Since $a\overline{\frac{1}{a}}=0$ in $Q/R$, it remains only to be shown that $\ann^R(\overline{\frac{1}{a}})\subseteq Ra$ and $\ann^{Q/R}(a)\subseteq R\overline{\frac{1}{a}}$.  Suppose that $\overline{\frac{1}{a}}b=0$.  Then $\frac{b}{a}\in R$, which shows that $b\in Ra$.  On the other hand, if $\overline{\frac{p}{q}}a=0$ in $Q/R$, then $\frac{pa}{q}\in R$.  Thus $\overline{\frac{p}{q}}=\frac{pa}{q}\overline{\frac{1}{a}}\in R\overline{\frac{1}{a}}$.
\end{proof}

\begin{lemma}\label{lem:comm2}
Let $R$ be a commutative B\'{e}zout ring with classical ring of quotients $Q$.  For any element $\overline{\frac{p}{q}}\in Q/R$, $\ann^R(\overline{\frac{p}{q}})=Rq$ and $\ann^{Q/R}(q)=R\overline{\frac{p}{q}}$.
\end{lemma}

\begin{proof}
We first claim that we may assume without loss of generality that $pR+qR=R$.  Since $R$ is B\'{e}zout, $pR+qR=bR$ for some $b\in R$.  Thus there are elements $s,t,x,y\in R$ such that $bs=p$, $bt=q$ and $px+qy=b$.  Thus $\frac{p}{q}=\frac{s}{t}$ and $sR+tR=R$ since $b$ is a nonzerodivisor.

Since it is immediate that $q\overline{\frac{p}{q}}=0$ in $Q/R$, it remains only to be shown that $\ann^R(\overline{\frac{p}{q}})\subseteq Rq$ and $\ann^{Q/R}(q)\subseteq R\overline{\frac{p}{q}}$.  Suppose that $a\overline{\frac{p}{q}}=0$.  Then $a\frac{p}{q}=r\in R$.  Since $pR+qR=R$, there are elements $u,v\in R$ such that $pu+qv=1$.  Thus $a=apu+aqv=qru+qav\in qR$, which implies that $\ann^R(\overline{\frac{p}{q}})\subseteq Rq$.  On the other hand, suppose that $\overline{\frac{c}{d}}q=0$.  By Lemma~\ref{lem:comm1}, $\overline{\frac{c}{d}}\in R\overline{\frac{1}{q}}$.  But $u\overline{\frac{p}{q}}=\overline{\frac{pu}{q}}=\overline{\frac{1-qv}{q}}=\overline{\frac{1}{q}}$.  Thus $\overline{\frac{c}{d}}\in R\overline{\frac{p}{q}}$, which shows that $\ann^{Q/R}(q)=R\overline{\frac{p}{q}}$.
\end{proof}

\begin{lemma}\label{lem:comm3}
Let $R$ be a commutative ring with classical ring of quotients $Q$.  If $a\in R$ is not a zerodivisor, then $(a,x)$ is morphic in $R\propto Q/R$ for every $x\in Q/R$.
\end{lemma}

\begin{proof}
Write $S=R\propto Q(R)/R$.  We will apply Lemma~\ref{lem:ann} to show that $\ann^S(a,x)=S(0,\overline{\frac{1}{a}})$ and $\ann^S(0,\overline{\frac{1}{a}})=S(a,x)$.  Since $a$ is not a zerodivisor, $\ann^R(a)=0$.  It is also clear that $(Q/R)a=Q/R$.  By Lemma~\ref{lem:comm1}, $\ann^R(\overline{\frac{1}{a}})=Ra$ and $\ann^{Q/R}(a)=R\overline{\frac{1}{a}}$.  Lemma~\ref{lem:ann} then gives the desired conclusion.
\end{proof}

\begin{lemma}\label{lem:comm4}
Let $R$ be a commutative B\'{e}zout ring with classical ring of quotients $Q$.  For any element $\overline{\frac{p}{q}}\in Q/R$, $(0,\overline{\frac{p}{q}})$ is morphic in $R\propto Q/R$.
\end{lemma}

\begin{proof}
Let $S=R\propto Q/R$.  We will show that $\ann^S(q,0)=S(0,\overline{\frac{p}{q}})$ and $\ann^S(0,\overline{\frac{p}{q}})=S(q,0)$.  As in the proof of Lemma~\ref{lem:comm3}, the result follows from Lemma~\ref{lem:ann} and Lemma~\ref{lem:comm2}.
\end{proof}

\begin{thm}\label{thm:red}
Suppose that $R$ is a commutative reduced B\'{e}zout weak Baer ring with classical ring of quotients $Q$.  Then $R\propto Q/R$ is morphic.
\end{thm}

\begin{proof}
Let $(a,x)\in R\propto Q/R$.  Since $R$ is a weak Baer ring, $\ann^R(a)=eR$ for some idempotent $e\in R$.  Write $f=1-e$, and let $S=eR$ and $T=fR$.  Then $R=S\times T$, and $R\propto Q/R=(S\propto Q(S)/S)\times (T\propto Q(T)/T)$.  Using this decomposition, $(a,x)=((0, xe),(af,xf)$.  By Lemma~\ref{lem:comm4}, $(0,xe)$ is morphic in $S\propto Q(S)/S$, and by Lemma~\ref{lem:comm3}, $(af,xf)$ is morphic in $T\propto Q(T)/T$ since $af$ is a nonzerodivisor in $Rf$.  Thus $(a,x)$ is morphic in $R\propto Q/R$.
\end{proof}

\begin{cor}\label{cor:moriffbaer}
Let $R$ be a commutative reduced ring with classical ring of quotients $Q$.  Then $R\propto Q/R$ is morphic if and only if $R$ is a B\'{e}zout weak Baer ring.
\end{cor}

\begin{proof}
Suppose that $R\propto Q/R$ is morphic.  By Proposition~\ref{prop:bez}, $R$ is B\'{e}zout.  Lemma~\ref{lem:prin} and Lemma~\ref{lem:baer} then show that $R$ is weak Baer.

On the other hand, if $R$ is B\'{e}zout and weak Baer, then Theorem~\ref{thm:red} implies that $R\propto Q/R$ is morphic.
\end{proof}

\begin{cor}\label{cor:bezdom}
If $R$ is a commutative B\'{e}zout domain with classical ring of quotients $Q$, then $R\propto Q/R$ is morphic.
\end{cor}

\begin{proof}
A commutative domain is both reduced and weak Baer. Corollary~\ref{cor:moriffbaer} then applies.
\end{proof}

In the case of an arbitrary domain, we can say more.
We begin with some general results about nonzerodivisors.  Recall that an element $a$ in a ring $R$ is called a left nonzerodivisor (respectively a right nonzerodivisor) if $\ann_r^R(a)=0$ (respectively $\ann_l^R(a)=0$).

\begin{prop}\label{prop:nzd}
Suppose that $R$ is a ring and that $M$ is a bimodule such that $R\propto M$ is left morphic.  If $a\in R$ is either a left or a right nonzerodivisor, then $Ma=M$ and there exists an element $m\in M$ such that $\ann_l^R(m)=Ra$ and $\ann_l^M(a)=Rm$.
\end{prop}

\begin{proof}
Write $S=R\propto M$, and consider the element $(a,0)\in S$.  Since $S$ is left morphic, there is an element $(b,m)\in S$ such that $\ann_l^S(a,0)=S(b,m)$ and $\ann_l^S(b,m)=S(a,0)$.  Then $ab=0=ba$.   Whether $a$ is a left or right nonzerodivisor, we can immediately conclude that $b=0$.  By Lemma~\ref{lem:ann}, we see that $Ma=M$ and that $\ann_l^R(m)=Ra$ and $\ann_l^M(a)=Rm$.
\end{proof}

\begin{cor}\label{cor:divring}
Let $R$ be a domain, and let $M$ be a bimodule. Suppose that there is an element $x\in M$ such that $\ann_l^R(x)=0$.  Then $R\propto M$ is left morphic and if and only if $R$ is a division ring and $R\propto M\cong R[t;\sigma ]/(t^2)$ for some ring endomorphism $\sigma$ of $R$.
\end{cor}

\begin{proof}
If $M$ contains an element $x$ such that $\ann_l^R(x)=0$, then Corollary~\ref{cor:tor} implies that $M\cong R(\sigma )$ for some ring endomorphism $\sigma$ of $R$, and Corollary~\ref{cor:ureg} shows that $R$ is unit regular. Thus $R$ is a division ring, and $R\propto M\cong R[t;\sigma ]/(t^2)$.

On the other hand, if $R$ is a division ring and $R\propto M\cong R[t;\sigma ]/(t^2)$ for some ring endomorphism $\sigma$, then $R\propto M$ is left morphic by \cite[Theorem 1]{cz:ext}.
\end{proof}

A module $M_R$ is called {\em divisible} (see \cite[p. 70]{lam:lmr}) if $x\in Ma$ for every pair of elements $a\in R$ and $x\in M$ such that $\ann_r^R(a)\subseteq\ann^R(x)$.  The next result characterizes morphic extensions over an arbitrary domain which is not a division ring (the case of a division ring having been covered previously).

\begin{thm}\label{thm:dom}Suppose that $R$ is a domain which is not a division ring and that $M$ is a bimodule.  Then $R\propto M$ is left morphic if and only if the following three conditions hold:
\begin{enumerate}
\item[{\rm (1)}] $M_R$ is divisible.
\item[{\rm (2)}] For every nonzero $a\in R$ there exists $m\in M$ such that $\ann_l^R(m)=Ra$ and $\ann_l^M(a)=Rm$.
\item[{\rm (3)}] For every $m\in M$ there exists a nonzero $a\in R$ such that $\ann_l^R(m)=Ra$ and $\ann_l^M(a)=Rm$.
\end{enumerate}
\end{thm}

\begin{proof}
Suppose that $R\propto M$ is left morphic.  Conditions (1) and (2) hold by Proposition~\ref{prop:nzd}.  By Corollary~\ref{lem:fund}, for every $m\in M$ there is an $a\in R$ as in condition (3), and $a$ is nonzero by Corollary~\ref{cor:divring}.

On the other hand, suppose that conditions (1), (2) and (3) hold.  Write $S=R\propto M$, and let $(a,x)\in S$.  If $a\ne 0$, then there is an element $m\in M$ such that $\ann_l^R(m)=Ra$ and $\ann_l^M(a)=Rm$ by condition (2).  Then $\ann_l^S(a,x)=0\propto\ann_l^M(a)=0\propto Rm=S(0,m)$, and $\ann_l^S(0,m)=\ann_l^R(m)\propto M=Ra\propto M=S(a,x)$ since $Ma=M$ by condition (1).

Suppose now that $a=0$.  By condition (3), there is a nonzero element $b\in R$ such that $\ann_l^R(x)=Rb$ and $\ann_l^M(b)=Rx$.  Since $b\ne 0$, then we see as above that $\ann_l^S(0,x)=S(b,0)$ and $\ann_l^S(b,0)=S(0,x)$.
\end{proof}

\begin{remark}
The conditions of Theorem~\ref{thm:dom} can be stated succinctly as follows.  If $R$ is a domain that is not a division ring and $M$ is a bimodule, then $R\propto M$ is left morphic if and only if for every nonzero $a\in R$ there is an element $m\in M$ and for every $m\in M$ there is a nonzero $a\in R$ such that the following sequence of left $R$-modules is exact:

$$0\rightarrow {}_RR\stackrel{\cdot a}{\rightarrow} {}_RR\stackrel{\cdot m}{\rightarrow} {}_RM\stackrel{\cdot a}{\rightarrow}{}_RM\rightarrow 0.$$
\end{remark}

The next result is a strengthening of Theorem~\ref{thm:str} in the domain case.

\begin{cor}\label{cor:lattice}
Suppose that $R$ is a domain which is not a division ring and that  $M$ is a bimodule such that $R\propto M$ is morphic.  The map $\mathcal{F}$, as defined in Theorem~\ref{thm:str}, is an inclusion reversing bijection between the set of cyclic right submodules of $M_R$ and the set of nonzero principal left ideals of $R$.
\end{cor}

\begin{proof}
Everything except the surjectivity of $\mathcal{F}$ follows from Theorem~\ref{thm:str}, suitably adapted to the case where ${}_RM$ is torsion.  The surjectivity of $\mathcal{F}$ follows from condition (2) of Theorem~\ref{thm:dom}.
\end{proof}

Given that Corollary~\ref{cor:lattice} provides so much control over the submodule lattice of $M$, we ask the following question.

\begin{question}\label{thm:convbezdom}
Suppose that $R$ is a commutative B\'{e}zout domain that is not a field and that $M$ is a bimodule such that $rm=mr$ for all $m\in M$ and all $r\in R$.  If $R\propto M$ is morphic, must $M\cong Q(R)/R$?
\end{question}

\begin{remark}In the situation above, we note that Corollary~\ref{cor:bezdom} and Corollary~\ref{cor:lattice} provide an isomorphism between the lattice of finitely generated submodules of $Q(R)/R$ and that of $M$.  We further remark that Theorem~\ref{thm:dom} implies that corresponding submodules are isomorphic.  Since a module is the direct limit of its finitely generated submodules, we conjecture that the question can be answered in the affirmative.
\end{remark}

It is further proved in \cite[Corollary 15]{cz:ext} that $\mathbb{Z}\propto\mathbb{Q}/\mathbb{Z}$ is {\em strongly} morphic.  Although we cannot yet extend this result to include an arbitrary commutative B\'{e}zout domain, we can extend it to the case of an elementary divisor domain.

A ring $R$ is called an elementary divisor ring if for every $m\times n$ matrix $A$ over $R$, there are invertible matrices $P$ and $Q$ such that $PAQ$ is diagonal with diagonal entries $d_i$ such that $d_i$ divides $d_{i+1}$ for every $i$, $1\le i<{\rm min}\lbrace m,n\rbrace$.  It is known that every elementary divisor domain is a B\'{e}zout domain, but the converse remains unsettled.  Therefore, although the next theorem holds only for elementary divisor domains, we are not able to appeal to a known example of a B\'{e}zout domain that is not an elementary divisor domain in order to determine whether the result is sharp.

\begin{thm}\label{thm:domstr}
If $R$ is an elementary divisor domain and $M$ is a bimodule such that $R\propto M$ is morphic, then $R\propto M$ is strongly morphic.  In particular, $R\propto Q/R$ is strongly morphic, where $Q$ denotes the classical ring of quotients of $R$.
\end{thm}

\begin{proof}
By Corollary~\ref{cor:strong}, the result is true if $R$ is a field.  For the remainder of the proof, we therefore assume that $R$ is not a field.

Let $B\in\mathbb{M}_n(R\propto M)$.  We claim that there are invertible matrices $U$ and $V$ such that $UBV$ is diagonal.  Then the matrix $UBV$ is morphic since it is a diagonal matrix all of whose entries are morphic.  By \cite[Lemma 3]{nisc:morp}, $B$ is morphic.

We now turn our attention to the claim.  Suppose that $B\in\mathbb{M}_n(R\propto M)$.  Exploiting the isomorphism between $\mathbb{M}_n(R\propto M)$ and $\mathbb{M}_n(R)\propto\mathbb{M}_n(M)$, we may write $B=(B_1,B_2)$.  Since $R$ is an elementary divisor domain, there are invertible matrices $P_1,Q_1\in\mathbb{M}_n(R)$ such that $P_1B_1Q_1=D_1$ is a diagonal matrix.  We may further assume that all of the nonzero diagonal entries precede those which are zero.  Writing $P=(P_1,0)$ and $Q=(Q_1,0)$, the matrix $B$ is equivalent to the matrix $PBQ=(D_1,C_2)$ for some matrix $C_2\in\mathbb{M}_n(M)$.

Before further reducing $B$ to a diagonal matrix, we consider one special case.  Let $N=(n_{ij})\in\mathbb{M}_n(M)$.  Since $R\propto M$ is morphic, $M$ is B\'{e}zout by Proposition~\ref{prop:bez}.  Therefore, the submodule of $M$ generated by the entries $n_{ij}$ is generated by the single element $n\in M$, and there are elements $r_{ij}\in R$ such that $n_{ij}=r_{ij}n$.  Since $R$ is an elementary divisor domain, there are invertible matrices $X,Y\in\mathbb{M}_n(R)$ such that $X(r_{ij})Y$ is diagonal.  Therefore $(0,N)$ is equivalent to the diagonal matrix $(X,0)(0,N)(Y,0)$.

We are now in position to prove that a given $B=(B_1,B_2)\in\mathbb{M}_n(R\propto M)$ is equivalent to a diagonal matrix.  By our first reduction, we may assume that
$$B_1=\begin{bmatrix}D_1&0\\0&0\end{bmatrix}$$
is block diagonal with $D_1$ an $k\times k$ diagonal matrix with no nonzero diagonal entries.  By our second reduction, we may then further assume that
$$B_2=\begin{bmatrix}\ast &\ast\\\ast & D_2\end{bmatrix}$$
is diagonal in the lower right-hand block.

If $B_1$ has no nonzero entries, then we are done.  Assume, then, that the (1,1) entry of $B_1$ is $d\ne 0$.  By induction, we may assume that every off-diagonal entry of $B$ is zero, except perhaps for the first row and column.  Writing $B=(b_{ij})$, suppose that $b_{11}=(d,m_{11})$ and $b_{1j}=(0,m_{1j})$.  Since $M$ is divisible by Theorem~\ref{thm:dom} and $d\ne 0$, there is an element $x_{1j}\in M$ such that $dx_{1j}=m_{1j}$.  If we subtract $(0,x_{1j})$ times the first column from the $j^{th}$ column, then the new $(1,j)$ entry will be $b_{1j}-b_{11}(0,x_{1j})=(0,m_{1j})-(d,m_{11})(0,x_{1j})=(0,0)$.  Continuing in this fashion, we can right and left multiply $B$ by elementary matrices to bring it to a diagonal form.  This proves the claim and the result.

\end{proof}

\begin{remark}
A result of Kaplansky (see \cite[p. 115]{fs:nonnoeth}) states that any diagonal matrix over a B\'{e}zout domain is equivalent to a diagonal matrix with $d_i$ dividing $d_{i+1}$ for all $i$, which brings us essentially back to the case of an elementary divisor domain.  Any possible extension of Theorem~\ref{thm:domstr} to the case of B\'{e}zout domains which may not be elementary divisor domains cannot, therefore, rely on being able to bring the matrices in question into a diagonal form.
\end{remark}

\bibliographystyle{alpha}
\bibliography{morphic}

\begin{thebibliography}{NSC04b}
\bibitem[Jon70]{jon:perf}
David Jonah.
\newblock Rings with the minimum condition for principal right ideals have the
  maximum condition for principal left ideals.
\newblock {\em Math. Z.}, 113:106--112, 1970.
\end{thebibliography}

\noindent
Department of Mathematics\\
Wellesley College\\
Wellesley, MA 02481\\
USA\\
Email: {\tt adiesl@wellesley.edu} \\

\noindent
Center for Communications Research\\
4320 Westerra Court\\
San Diego, CA 92126-1967\\
USA\\
Email: {\tt dorsey@ccrwest.org} \\

\noindent
Department of Mathematics and Statistics \\
Bowling Green State University \\
Bowling Green, OH 43403 \\
USA \\
Email: {\tt warrenb@bgsu.edu}

\end{document}